\documentclass[10pt]{article}
\usepackage{IEEEtrantools}
\usepackage{mathtools, nccmath}
\usepackage{latexsym}
\usepackage{amsfonts}
\usepackage{amssymb}
\usepackage{graphicx}
\usepackage{amsmath,amsfonts,amsthm}
\usepackage{mathrsfs}
\usepackage{enumerate}
\usepackage{dsfont}
\usepackage{hyperref}
\usepackage{pst-node}
\usepackage{cleveref}
\usepackage{tikz-cd}

\usepackage{multirow}
\usepackage{tikz}
\usetikzlibrary{arrows,decorations.markings}
\usepackage{xcolor}
\usepackage{bookmark}
\usepackage{hyperref}
\hypersetup{
     colorlinks   = true,
     citecolor    = blue
}

\newcommand\restr[2]{{
  \left.\kern-\nulldelimiterspace 
  #1 
  \littletaller 
  \right|_{#2} 
  }}

\newcommand{\littletaller}{\mathchoice{\vphantom{\big|}}{}{}{}}

\newcommand{\be}{\begin{equation}}
\newcommand{\ee}{\end{equation}}
\newcommand{\bea}{\begin{eqnarray}}
\newcommand{\eea}{\end{eqnarray}}
\newcommand{\bean}{\begin{eqnarray*}}
\newcommand{\eean}{\end{eqnarray*}}
\newcommand{\brray}{\begin{array}}
\newcommand{\erray}{\end{array}}
\newcommand{\biearray}{\begin{IEEEarray}{rCl}}
\newcommand{\eiearray}{\end{IEEEarray}}

\newtheorem{dfn}{Definition}[section]
\newtheorem{thm}[dfn]{Theorem}
\newtheorem{lmma}[dfn]{Lemma}
\newtheorem{ppsn}[dfn]{Proposition}
\newtheorem{crlre}[dfn]{Corollary}
\newtheorem{xmpl}[dfn]{Example}
\newtheorem{rmrk}[dfn]{Remark}

\newcommand{\bdfn}{\begin{dfn}\rm}
\newcommand{\bthm}{\begin{thm}}
\newcommand{\blmma}{\begin{lmma}}
\newcommand{\bppsn}{\begin{ppsn}}
\newcommand{\bcrlre}{\begin{crlre}}
\newcommand{\bxmpl}{\begin{xmpl}}
\newcommand{\brmrk}{\begin{rmrk}\rm}

\newcommand{\edfn}{\end{dfn}}
\newcommand{\ethm}{\end{thm}}
\newcommand{\elmma}{\end{lmma}}
\newcommand{\eppsn}{\end{ppsn}}
\newcommand{\ecrlre}{\end{crlre}}
\newcommand{\exmpl}{\end{xmpl}}
\newcommand{\ermrk}{\end{rmrk}}

\newcommand{\tr}{\mathrm{tr}}

\def \qed { \mbox{}\hfill
$\Box$\vspace{1ex}}

\title{Sections and Chapters}

\begin{document}

\tikzset{->-/.style={decoration={
  markings,
  mark=at position #1 with {\arrow{>}}},postaction={decorate}}}
  \tikzset{-<-/.style={decoration={
  markings,
  mark=at position #1 with {\arrow{<}}},postaction={decorate}}}

\author{\sc{Keshab Chandra Bakshi,\,\,Satyajit Guin,\,\,Biplab Pal}}

\title{Von Neumann entropy of the angle operator between a pair of intermediate subalgebras}
\maketitle

%%%%%%%%%%%%%%%%%%%%%%%%%%%%%%%%%%
%%%%%%%%%%%   ABSTRACT    %%%%%%%%%%%%%%%%
%%%%%%%%%%%%%%%%%%%%%%%%%%%%%%%%%%

\begin{abstract}
Given a pair of intermediate $C^*$-subalgebras of a unital inclusion of simple $C^*$-algebras with a conditional expectation of finite Watatani index, we discuss the corresponding angle operator and its Fourier transform. We  provide a calculable formula for the von Neumann entropy of the (Fourier) dual angle operator for a large class of quadruple of simple $C^*$-algebras. \end{abstract}
\bigskip

{\bf AMS Subject Classification No.:} {\large 46}L{\large 37}\,, {\large 47}L{\large 40}\,, {\large 46}L{\large 05}\,, {\large 43}A{\large 30}\,, {\large{94}A{\large 15}\,.
\smallskip

{\bf Keywords.} simple $C^*$-algebra, Watatani index, Fourier transform, angle operator, von Neumann entropy

\hypersetup{linkcolor=blue}
\bigskip

\section{Introduction}

In recent years, the study of the symmetries of an inclusion of $C^*$-algebras becomes an active area of research. In this short article, we focus on a unital inclusion of simple $C^*$-algebras $B\subset A$ with a conditional expectation of finite \textit{Watatani index} (a generalization of Jones index) and a pair of intermediate $C^*$-subalgebras $B\subset C, D\subset A$. Motivated by  \cite{SW}, we may consider the corresponding \textit{angle operator} $\Theta$ between $C$ and $D$, and using the \textit{Fourier transform} on the \textit{relative commutants} of $B\subset A$, as developed in \cite{BG}, we may consider the Fourier dual $\mathcal{F}(\Theta)$ of the angle operator. Our goal is to compute the \textit{von Neumann entropy} of the angle operator and its Fourier dual. We observe that the subalgebras $C$ and $D$ are `orthogonal' (in others words, forms a `commuting square') if and only if the von Neumann entropy of the angle operator vanishes. The main result of this article is the following theorem (see \Cref{main}).
\smallskip
 
\noindent \textbf{Theorem A:} Let $B\subset A$ be an irreducible inclusion of simple unital $C^*$-algebras with a conditional expectation of index-finite type and suppose that $C,D$ are two intermediate unital $C^*$-subalgebras. Then,
\begin{center}
$H\big(|\mathcal{F}(\Theta)|^2\big)=\frac{2}{\sqrt{[A:B]_0}} \eta\big(\delta \tr(e_Ce_D)\big).$
\end{center}
\smallskip
 
We remark that generalizing the above formula beyond the irreducible situation seems difficult at the moment. However, in the case of `co-commuting squares' (a dual notion of commuting square), we have explicit formula (see \Cref{cocom}).
\smallskip

\noindent \textbf{Theorem B:} Let $(B\subset C,D\subset A)$ be a quadruple of simple unital $C^*$-algebras with a conditional expectation from $A$ onto $B$ of index-finite type. If the quadruple is a co-commuting square, then
\begin{center}
$H\big(|\mathcal{F}(\Theta)|^2\big)=\frac{2}{\sqrt{[A:B]_0}} \hspace{0.1 cm} \eta\bigg(\displaystyle \frac{\sqrt{[A:B]_0}}{[A:C]_0[A:D]_0}\bigg).$
\end{center}

 %%%%%%%%%%%%%%%%%%%%%%%%%%%%%%%%%%%%%
 %%%%%%%%%%%%%%%%%%%%%%%%%%%%%%%%%%%%%
 %%%%%%%%%%%%%%%%%%%%%%%%%%%%%%%%%%%%%
 
\section{Preliminaries}\label{plms}

We first recall Watatani's  $C^*$-index theory for a unital inclusion of simple $C^*$-algebras following \cite{W}, and then briefly touch upon the Fourier theory developed in \cite{BG} for the relative commutants of that inclusion. Consider a unital inclusion of $C^{*}$-algebras $B\subset A$ and a conditional expectation $E:A \rightarrow B$. $E$ is said to be of  index finite type if there exist a quasi-basis  $ \{ \lambda_1 , \lambda_2 ,\cdots , \lambda_n\} \subset A $. In this case the Watatani index of $E$ is given by $\text{Ind}_w(E)=\sum\limits_{i=1}^{n} \lambda_i\lambda^{*}_i .$ For an inclusion $B\subset A$ of simple $C^*$-algebras with a contional expectation of index-finite type, there always exists a minimal conditional expectation $E^A_B$. The minimal index of $B\subset A$, denoted by $[A:B]_0$, is defined as $\text{Ind}_w(E^A_B)$. In the sequel, we always deal with an inclusion of simple unital $C^*$-algebras. We denote the $C^*$-basic construction of the inclusion $B\subset A$ by $A_1$. Iterating the $C^*$-basic construction, we get a tower of simple unital $C^*$-algebras 
$B\subset A \subset A_1 \subset A_2 \subset \cdots \subset A_k \subset \cdots$ and obtain the unique minimal condition expectation $E^{A_k}_{A_{k-1}} : A_k \rightarrow A_{k-1}$ for each $k \geq 0$. Let $e_{k+1}$ denote the Jones projection that implements the basic construction of the inclusion $A_{k-1} \subset A_k$ with respect to $E^{A_k}_{A_{k-1}}$. The relative commutants $B'\cap A_k$ are finite-dimensional. For each $k\geq 0$, $E^{A}_B \circ E^{A_1}_A \circ \cdots \circ E^{A_k}_{A_{k-1}} \big|_{B'\cap A_k}$ is a faithful tracial state on $B'\cap A_k$, to be denoted by $\tr_k$. We often drop `$k$' and denote $\tr_k$ simply by $\tr$. The unique trace preserving conditional expectation $E^{B'\cap A_{k}}_{A'\cap A_{k}} : B'\cap A_{k} \rightarrow A'\cap A_{k}$ is  given by $$E^{B'\cap A_{k}}_{A'\cap A_{k}}\big(x\big)=\frac{1}{[A:B]_0}\sum\limits_{i} \lambda_i x\lambda^{*}_i\,.$$ Following \cite{BG} (see also \cite{BGS}), we quickly recall the notion of `Fourier transform' and `convolution' product on the relative commutants. For each $k \geq 0$, the {\em Fourier transform} $\mathcal{F}_k : B'\cap A_k \rightarrow A'\cap A_{k+1}$ is defined as 
$$\mathcal{F}_k(x)=\delta^{k+2} E^{B'\cap A_{k+1}}_{A'\cap A_{k+1}} \big(xe_{k+1}e_k \cdots e_2 e_1) \text{  for all } x\in B'\cap A_k;$$ and the {\em inverse Fourier transform} $\mathcal{F}^{-1}_k :  A'\cap A_{k+1} \rightarrow B'\cap A_k $ is defined as $$\mathcal{F}^{-1}_k(x) =\delta^{k+2} E^{A_{k+1}}_{A_k}\big(xe_1e_2 \cdots e_k e_{k+1}\big) \text{  for all  } x\in A'\cap A_{k+1} .$$
Here we use the term ``inverse" in the sense that $\mathcal{F}_k\circ\mathcal{F}^{-1}_k=\text{id}_{A'\cap A_{k+1}}$and $\mathcal{F}^{-1}_k\circ\mathcal{F}_k=\text{id}_{B'\cap A_{k}}$. For simplicity, we shall use $\mathcal{F}$ for $\mathcal{F}_1$.

\begin{thm}[\cite{BG}, Theorem 3.5]\label{th5}
$\mathcal{F}$ and $\mathcal{F}^{-1}$ are isometries with respect to the norm defined by $||x||_2=\tr\big(x^{*}x\big)$.
\end{thm}

For $x,y\in B'\cap A_1$, the {\em convolution} product of $x$ and $y$, to be denoted by $x \star y$, is defined by the following:
\begin{equation*}
x\star y:=\mathcal{F}^{-1}\big(\mathcal{F}(y)\mathcal{F}(x)\big).
\end{equation*}
Similarly, for $w,z \in A'\cap A_2$, their {\em convolution} product is defined as
\begin{equation*}
w\star z:=\mathcal{F}\big(\mathcal{F}^{-1}(z)\mathcal{F}^{-1}(w)\big).
\end{equation*}
The convolution product is associative (\cite{BG}, Lemma 3.20). For $x,y\in B'\cap A_1$, by Proposition $3.8$ in \cite{BGS}, we have
\begin{equation*}
( x\star y)^{*} = x^{*} \star y^{*};
\end{equation*}
and similarly for $w,z\in A'\cap A_2$, we have
\begin{equation}\label{eqn4}
(w\star z)^{*} = w^{*} \star z^{*}.
\end{equation}
Suppose that $C$ is an intermediate simple unital $C^*$-subalgebra of $B\subset A$. By \cite{I}, there exist minimal conditional expectations $E^{C}_B$ from $C$ onto $B$ and $E^{A}_C$ from $A$ onto $C$ such that $E^{A}_B = E^{C}_B \circ E^{A}_C$. Indeed, $E^C_B$ is given by restricting $E^A_B$ onto $C$. Let $C_1$ denote the $C^*$-basic construction of the inclusion $C\subset A$ with Jones projection $e_C$ corresponding to the minimal conditional expectation $E^{A}_C$. Below we list some useful results, the proofs of which follow along the same line of argument as in \cite{BG}, Section 4.

\blmma[\cite{BG}, Lemma 4.2]\label{lma4}
Let $B\subset C \subset A$ be as discussed above. Then, we have the following:
\begin{enumerate}
\item $C_1\subset A_1$ is an inclusion of simple unital $C^*$-algebras with common identity. The dual conditional expectations $E^{A_1}_A$ and $E^{C_1}_A$ are minimal.  $E^{A_1}_A$ and $E^{C_1}_A$ must satisfy $E^{A_1}_A = E^{C_1}_A E^{A_1}_{C_1}$, and hence $E^{A_1}_A\big|_{C_1} = E^{C_1}_A$.
\item The unique $\tr$-preserving conditional expectation from $B'\cap A_1$ onto $B'\cap C_1$ is given by $E^{A_1}_{C_1}\big|_{B'\cap A_1}$.
\item The tracial state on the relative commutant $C^{\,\prime}\cap C_1$ induced by the inclusion $C\subset A \subset C_1$ is the restriction of the tracial state on $B'\cap A_1$ induced by the inclusion $B\subset A \subset A_1$.
\end{enumerate}
\elmma

\bppsn[\cite{BG}, Lemma 4.4(2) and Proposition 4.6]\label{ppn9} \label{fact2}
Let $B\subset C \subset A$ be as in Lemma \ref{lma4}. Then, we have the following: 
\begin{enumerate}
\item $[A_1:C_1]_0 = [C:B]_0$;
\item $E^{A_1}_{C_1}(e_1)=\frac{1}{[C:B]_0}e_C$;
\item $\mathcal{F}(e_C) =\frac{\sqrt{[A:B]_0}}{[A:C]_0} e_{C_1}$ and $\mathcal{F}(e_1) = \frac{1}{\sqrt{[A:B]_0}}$.
\end{enumerate}
\eppsn
\smallskip
  
In what follows, we shall fix a couple of notations  that will be used in the sequel.
\smallskip

\noindent\textbf{Notation:} 
\begin{enumerate} 
\item The minimal index $[A:B]_0$ will be denoted by $\delta^2$. Also, $r=\frac{[C:B]_0}{[A:D]_0}=\frac{[D:B]_0}{[A:C]_0},\,\tau=\frac{1}{[A:B]_0}$, and $\tau_C =\frac{1}{[A:C]_0}$.
\item $\kappa^{+}_0 =$min $\{\tr(p) : p\in \mathcal{P}(B'\cap A)\}$, $\kappa^{-}_0 =$min $\{\tr(q) : q\in \mathcal{P}(A'\cap A_1)\}$, and $\kappa_0=\sqrt{\kappa^{+}_0\kappa^{-}_0}$.
\item The shorthand notation $(B\subset C,D \subset A)$ denotes the following quadruple
\[\begin{matrix}
D &\subset & A\\
\cup & &\cup\\
B &\subset & C
\end{matrix}\]
of $C^*$-algebras.
 \end{enumerate}
 
%%%%%%%%%%%%%%%%%%%%%%%%%%%%%%%%%%%%%
%%%%%%%%%%%%%%%%%%%%%%%%%%%%%%%%%%%%%
%%%%%%%%%%%%%%%%%%%%%%%%%%%%%%%%%%%%%
 
\section{Fourier transform of the angle operator and its von Neumann entropy} \label{entropy}

Throughout this section, $(B\subset C,D \subset A)$ denotes a quadruple of simple unital $C^*$-algebras with a conditional expectation from $A$ onto $B$ of index-finite type. Given the quadruple ($B\subset C,D \subset A$), there exist unique mininimal conditional expectations $E^{A}_C$ and $E^{A}_D$. Let $e_C$ and $e_D$ be the Jones projection corresponding to  $E^{A}_C$ and $E^{A}_D$ respectively. By Lemma \ref{lma4}, the new quadruple ($A\subset C_1,D_1 \subset A_1$) again forms a quadruple of simple unital $C^*$-algebras, where $B_1 ,C_1 $ and $D_1$ are $C^*$-basic constructions corresponding to $E^{A}_B$, $E^{A}_C$ and $E^{A}_D$ respectively. The quadruple ($B\subset C,D \subset A$) will be called {\em commuting square} if $E^{A}_CE^{A}_D=E^{A}_DE^{A}_C=E^{A}_B $. The same quadruple will be called {\em co-commuting} if the dual quadruple ($A\subset C_1,D_1 \subset A_1$) is a commuting square.
  \smallskip
  
  Motivated by \cite{SW, LX}, we define the angle operator between a pair of intermediate simple $C^*$-subalgebras as follows.
\bdfn\label{dfn3}
Let $B \subset A$ be an inclusion of simple unital $C^*$-algebras, and $C,D$ be two intermediate (unital) simple $C^{*}$-algebras. The angle operator, denoted by $\Theta$, is defined as follows: $$\Theta = e_C e_D.$$ 
\edfn

\begin{dfn}[\cite{JZW}]\label{dfn1}
For $x\in B' \cap A_1$, the von Neumann entropy of $|x|^{2}$ is defined as the following quantity: $$H(|x|^{2})=tr(\eta
(|x|^2)) .$$
\end{dfn}

\bppsn\label{th7}
Suppose that $(B\subset C,D \subset A)$ is a quadruple of simple unital $C^*$-algebras with a conditional expectation from $A$ onto $B$ of index-finite type. Then,
\begin{center}
$H(|\Theta|^2)+ H(|\mathcal{F}(\Theta)|^2)\geq\frac{2{\kappa}_0}{\delta}\eta\big(\frac{\delta}{{\kappa}_o} \tr(e_Ce_D)\big).$
\end{center}
\eppsn
\begin{prf}
Use Theorem $4.10$ in \cite{BGS}.\qed
\end{prf}

\begin{ppsn}
If $(B\subset C,D\subset A)$ is a commuting square of simple unital $C^*$-algebras with a conditional expectation from $A$ onto $B$ of index-finite type, then  $H(|{\Theta}|^2)=0.$ Furthermore, in this case $H\big(|\mathcal{F}(\Theta)|^2\big)=\eta\big([A:B]^{-1}_0\big).$
\end{ppsn}
\begin{prf}
Since $e_Ce_D=e_De_C=e_1$, we get that $\lvert {\Theta}\rvert^2= e_1$. Therefore, $\eta(|{\Theta}|^2)=0$ and so, $H(\lvert{\Theta}\rvert^2)= 0.$
Also, as  $\mathcal{F}(e_1)=\frac{1}{\sqrt{[A:B]_0}}$, we have $\lvert \mathcal{F}(\Theta)\rvert ^2= \frac{1}{[A:B]_0}.$\qed
\end{prf}

The converse of the above theorem is also true.

\begin{ppsn}
Let $(B\subset C,D\subset  A)$ be an irreducible inclusion of simple unital $C^*$-algebras with a conditional expectation $E$ from $A$ onto $B$ of index-finite type. If $H(|{\Theta}|^2)=0$, then $$e_Ce_D = e_De_C = e_{C\cap D}.$$
\end{ppsn}
\begin{prf}
Observe that $0\leq e_De_Ce_D \leq 1 \text{ and } \eta \text{ is positive in}~ [0,1]$. Now, $H(|{\Theta}|^2)=0$ implies that $\eta(e_De_Ce_D)=0$, and by (\cite{PP}, Page $74$) we conclude that $e_De_Ce_D$ is a projection in $ B^{\prime}\cap A_1.$ Since $B'\cap A_1$ is finite-dimensional, $(e_De_Ce_D)^n$ converges to $e_C \wedge e_D$ in the norm topology. Also since $e_C\wedge e_D = e_{C\cap D}$ (\cite{BGJ}, Proposition $4.1$), we immediately get that $e_De_Ce_D=e_{C\cap D}$, and hence $e_De_C =e_C e_D= e_{C\cap D}$.\qed
\end{prf}

To provide a formula for the von Neumann entropy of the dual angle operator beyond the commuting square situation, let us recall the auxiliary operators associated to a quadruple of simple $C^*$-algebras, along with various properties of them, as developed in \cite{BG}.
\smallskip

Suppose that $\{ {\gamma}_i : 1 \leq i \leq m\}$ and $ \{ {\delta}_j : 1 \leq j \leq n\}$ are quasi-bases of $E^{C}_B$ and $E^{D}_B$ respectively. We can define two auxiliary operators $p(C,D)\in C'\cap D_1$ and $q(C,D)\in D'\cap C_1$ by
\begin{IEEEeqnarray}{lCl}
p(C,D):=\sum\limits_{i,j}\gamma_i \delta_j e_1 {\delta}^{*}_j {\gamma}^{*}_i\qquad\mbox{ and }\qquad q(C,D):=\sum\limits_{i,j} \delta_j \gamma_i e_1 {\gamma}^{*}_i {\delta}^{*}_j\,.
\end{IEEEeqnarray}
 
\blmma[\cite{BG}, Proposition $5.15$]\label{ppn10}
Let $B\subset A$ be an irreducible inclusion of simple unital $C^*$-algebras with a conditional expectation from $A$ onto $B$ of index-finite type. If $C$ and $D$ are two intermediate unital $C^*$-subalgebras of $B\subset A$, then $\frac{1}{t} p(C,D)$ and $\frac{1}{t} q(C,D)$ both are projections, where $t=[A:B]_0\tr\big(e_Ce_D\big)$.
\elmma

\blmma\label{ppn4}
If $(B\subset C,D\subset A)$ is a commuting square of simple  unital $C^*$-algebras with a conditional expectation from $A$ onto $B$ of index-finite type,  then $p(C,D)$ and $q(C,D)$ both are projections.
\elmma
\begin{prf} 
The proof follows along the same line of argument as in Proposition $2.20$ in \cite{BDLR} and hence we omit it.\qed
\end{prf}

The following result was proved in \cite{BG}  in the irreducible case only.

\blmma\label{ppn3}
If $(B\subset C,D \subset A)$ is a quadruple of simple unital $C^*$-algebras with a conditional expectation from $A$ onto $B$  of index-finite type, then $$ p(C,D)= [D:B]_0 \hspace{0.1 CM}E^{A_1}_{D_1}(e_C)\,,$$ and $$ q(C,D)= [C:B]_0 \hspace{0.1CM}E^{A_1}_{C_1}(e_D)\,.$$
\elmma
\begin{prf}
By \cite{BG}[Proposition $4.2(4)$], we know that $E^{A_1}_{C_1}\big|_{B'\cap A_1}$ is the unique trace preserving conditional expectation from $B'\cap A_1$ onto $B'\cap C_1$. For any $w=\sum\limits_{i} x_i e_C y_i \in B'\cap C_1,$ where $x_i , y_i \in A$ for all $i$, we have
\begin{IEEEeqnarray*}{lCl}
\tr(e_D w) &=&\tr\big(e_D\big(\textstyle\sum_ix_i e_C y_i\big)\big)\\
&=&\tr\big(\big(\textstyle\sum\limits_{k} \delta_k e_1 \delta^{*}_k\big)\big(\textstyle\sum\limits_{i} x_i e_C y_i\big)\big)\hspace{3.4cm} (\text{by \cite{BG}, Lemma $5.11$})\\
&=&\tr\big(\textstyle\sum\limits_{k,i} \delta_k e_1 e_C\delta^{*}_k x_i e_C y_i\big)\\
&=&E^{A}_{B} E^{C_1}_{A} E^{A_1}_{C_1}\big(\textstyle\sum\limits_{k,i} \delta_k e_1  E^{A}_{C}\big(\delta^{*}_k x_i\big)  y_i\big)\\
&=&E^{A}_{B} E^{C_1}_{A} \big(\textstyle\sum\limits_{k,i} \delta_k E^{A_1}_{C_1}\big(e_1\big)  E^{A}_{C}\big(\delta^{*}_k x_i\big)  y_i\big)\\
&=&\frac{1}{[C:B]_0}E^{A}_{B} E^{C_1}_{A} \big(\textstyle\sum\limits_{k,i} \delta_k  e_C E^{A}_{C}\big(\delta^{*}_k x_i\big)  y_i\big)\hspace{1.8cm} (\text{by \cite{BG}, Lemma $4.4(2)$})\\
&=&\frac{1}{[C:B]_0} E^{A}_{B}  \big(\textstyle\sum\limits_{k,i} \delta_k  E^{C_1}_{A}\big(e_C\big) E^{A}_{C}\big(\delta^{*}_k x_i\big)  y_i\big)\\
&=&\frac{1}{[C:B]_0[A:C]_0} E^{A}_{B}\big(\textstyle\sum\limits_{k,i} \delta_k   E^{A}_{C}\big(\delta^{*}_k x_i\big)  y_i\big)\,.\hspace{1.4 cm}(\mbox{since  }E^{C_1}_{A}\big(e_C\big)=[A:C]_0^{-1})\end{IEEEeqnarray*}
On the other hand, we have the following       
\begin{IEEEeqnarray*}{lCl}
\tr(q(C,D) w) &=& \tr\big(\big(\textstyle\sum\limits_{k} \delta_k e_C \delta^{*}_k\big)\big(\textstyle\sum\limits_{i} x_i e_C y_i\big)\big)\hspace{1.9cm} (\mbox{by \cite{BG}, Lemma $5.11$})\\
&=& \tr\big(\textstyle\sum\limits_{k,i} \delta_k e_C \delta^{*}_k  x_i e_C y_i\big)\\
&=&E^{A}_{B} E^{A_1}_{A} \big(\textstyle\sum\limits_{k,i} \delta_k E^{A}_{C}\big( \delta^{*}_k  x_i \big) e_C y_i\big)\\
&=&E^{A}_{B} \big(\textstyle\sum\limits_{k,i} \delta_k E^{A}_{C}\big( \delta^{*}_k  x_i \big) E^{A_1}_{A}\big( e_C\big) y_i\big)\\
&=&\frac{1}{[A:C]_0}E^{A}_{B} \big(\textstyle\sum\limits_{k,i} \delta_k E^{A}_{C}\big( \delta^{*}_k  x_i \big)  y_i\big). \hspace{2 cm}\\
\end{IEEEeqnarray*}
Therefore, we get
\[
\tr(q(C,D) w)=[C:B]_0 \tr(e_D w),
\]
and so
\[
q(C,D)=[C:B]_0\,E^{A_1}_{C_1}(e_D).
\]  
The proof for $p(C,D)$ is similar.\qed
\end{prf}

\begin{crlre}\label{fact1}
If $(B\subset C,D \subset A)$ is a quadruple of simple unital $C^*$-algebras with a conditional expectation from $A$ onto $B$ of index-finite type, then $$ \tr \big(p(C,D)\big)= \tr\big(q(C,D)\big) = r.$$
\end{crlre}
\begin{prf}
By \Cref{ppn3}, we have $q(C,D)= [C:B]_0 \hspace{0.1CM}E^{A_1}_{C_1}(e_D).$ Thus,
\begin{IEEEeqnarray*}{lCl}
\tr\big(q(C,D)\big)&=& [C:B]_0\,\tr\big(E^{A_1}_{C_1}(e_D)\big)\\
&=& [C:B]_0\,\tr\big(e_D\big)=r.
\end{IEEEeqnarray*}
The proof for $p(C,D)$ is similar.\qed
\end{prf}

\bppsn\label{ppn1}
If $(B\subset C,D \subset A)$ is a quadruple of simple unital $C^*$-algebras with a conditional expectation from $A$ onto $B$ of index-finite type, then $$p(C ,D) = \delta \hspace{0.1cm}e_C \star e_D$$ and
$$q(C ,D) = \delta \hspace{0.1cm}e_D \star e_C .$$
\eppsn
\begin{prf}
By \cite{BGS}[Lemma 4.15], we know that for any $x,y\in B^{\prime}\cap A_1$, we have
\[
E^{A_2}_{A_1}\big(\mathcal{F}(x) y\mathcal{F}(x)^*\big)=\frac{1}{\delta} y\star xx^*.
\]
Putting $x=e_D$ and $y=e_C$, we immediately obtain $\delta e_C\star e_D={[D:B]^2_0} E^{A_2}_{A_1}(e_{D_1}e_Ce_{D_1})$, thanks to \Cref{ppn9}. Since $e_{D_1}e_Ce_{D_1}=E^{A_1}_{D_1}(e_C)e_{D_1}$ and $E^{A_2}_{A_1}(e_{D_1})={[A_1:D_1]}^{-1}_0={[D:B]}^{-1}_0$ (see \Cref{fact2}), we observe that  $\delta e_C\star e_D= {[D:B]_0} E^{A_1}_{D_1} (e_C)$. The proof is complete once we apply \Cref{ppn3}. The proof for $q(C,D)$ is obtained by interchanging $C$ and $D$.\qed
\end{prf}
\blmma\label{ppn2}
If $(B\subset C,D \subset A)$ is a quadruple of simple unital $C^*$-algebras with a conditional expectation from $A$ onto $B$ of index-finite type, then
\begin{center}
$\tr(e_{C_{1}}e_{D_{1}})=\frac{[A:C]_0}{[D:B]_0}\hspace{0.2cm} \tr(e_Ce_D).$
\end{center} 
\elmma
\begin{prf}
Using \Cref{ppn9} and by the multiplicativity of the minimal index, we observe that the following equalities hold true:
\begin{IEEEeqnarray*}{lCl}
\tr(e_{C_{1}}e_{D_{1}}) &=& \big<e_{C_{1}} , e_{D_{1}} \big>\\
&=& \frac{{[A:C]_0}{[A:D]_0}}{[A:B]_0} \big< \mathcal{F}( e_C) ,\mathcal{F}( e_D) \big>\\
&=& \frac{[A:C]_0}{[D:B]_0}\big< e_C , e_D \big> \hspace{4.7cm} \text{(by \cite{BG}, Theorem 3.5)}\\
&=&\frac{[A:C]_0}{[D:B]_0}\hspace{0.2cm} \tr(e_Ce_D),
\end{IEEEeqnarray*}
which completes the proof.\qed
\end{prf}

\bthm\label{th1}\label{main}
Let $B\subset A$ be an irreducible inclusion of simple unital $C^*$-algebras with a conditional expectation from $A$ onto $B$ of index-finite type and suppose that $C,D$ are two intermediate (unital) $C^*$-subalgebras. Then,
\begin{center}
$H\big(|\mathcal{F}(\Theta)|^2\big)=\frac{2}{\delta} \eta\big(\delta \tr(e_Ce_D)\big).$
\end{center}
\ethm
\begin{prf}
First note that $  \mathcal{F}(\Theta)= r e_{D_1} \star e_{C_1}$, and therefore
\begin{IEEEeqnarray*}{lCl}
H\big(|\mathcal{F}(\Theta)|^2\big) &=& H\big({\mathcal{F}(\Theta)}^{*}{\mathcal{F}(\Theta)}\big)\\
&=& H\big((r e_{D_1}\star e_{C_1})^{*}(r e_{D_1}\star e_{C_1})\big)\\
&=& H\big(r^{2} ( {e^{*}_{D_1}} \star {e^{*}_{C_1}})(e_{D_1}\star e_{C_1})\big) \hspace{1cm}\mbox{(by \Cref{eqn4})}\\
&=& H\big(r^{2}  (e_{D_1}\star e_{C_1})^{2} \big).
\end{IEEEeqnarray*}
By \Cref{ppn1}, we have $q_1:=q(C_1 ,D_1) = \delta \hspace{0.1cm}e_{D_1} \star e_{C_1}$.
Since $B\subset A$ is an irreducible inclusion, it follows that the inclusion $A\subset A_1$ is also irreducible (see \cite{BG}[Proposition 3.2], for instance). Now, using \cite{BG}[Proposition 5.15] we know that $q^{2}_1= [A_1:A]_0 \tr ( e_{C_1}e_{D_{1}}) q_1$. Thus, we have
\begin{IEEEeqnarray*}{lCl}
(\delta \hspace{0.1cm}e_{D_1} \star e_{C_1})^{2} &=& [A:B]_0\delta \hspace{0.1cm}\tr ( e_{C_1}e_{D_{1}}) e_{D_1} \star e_{C_1}\\
&=& \frac{[A:C]_0{[A:B]^{\frac{3}{2}}_0}}{[D:B]_0}\hspace{0.2cm} \tr(e_Ce_D)e_{D_1} \star e_{C_1}.\hspace{1.5 cm}(\mbox{by \Cref{ppn2}})
\end{IEEEeqnarray*}
In other words, we see that $(e_{D_1} \star e_{C_1})^{2} =\displaystyle \frac{[A:C]_0\sqrt{[A:B]_0}}{[D:B]_0}\hspace{0.2cm} \tr(e_Ce_D)e_{D_1} \star e_{C_1}.$ If we put $f= \displaystyle \frac{\delta e_{D_1} \star e_{C_1}}{[A:C]_0[A:D]_0\tr(e_Ce_D)}$, then a straightforward calculation yields 
\begin{IEEEeqnarray}{lCl}\label{imp}
r^{2}  (e_{D_1}\star e_{C_1})^{2} = (\delta\,\tr(e_Ce_D))^{2}f.
\end{IEEEeqnarray}
Using \Cref{ppn2}, we can easily see that
\begin{center}
$f=\frac{\delta\, e_{D_1}\star\, e_{C_1}}{[A_1:A]_0\tr ( e_{C_1}e_{D_{1}})}.$
\end{center}
Thus, by \cite{BG}[Proposition 5.15] and Proposition \ref{ppn1}, we conclude that $f$ is a projection. Furthermore,
\begin{IEEEeqnarray*}{lCl}
\tr(f) &=&\displaystyle \frac{[C_1:A]_0}{[A_1:D_1]_0[A:C]_0[A:D]_0\tr(e_Ce_D)}\hspace{1cm}(\mbox{by \Cref{fact1}})\\
&=&\frac{1}{[A:B]_0\tr(e_Ce_D)} \hspace{4 cm}(\mbox{by \Cref{fact2}}).
\end{IEEEeqnarray*}
Therefore, using \Cref{imp}, we easily obtain the following equalities
\begin{IEEEeqnarray*}{lCl}
H\big(|\mathcal{F}(\Theta)|^2\big) &=& H\big(r^{2}\big(e_{D_1}\star e_{C_1}\big)^{2} \big)\\
&=&H\big(\big(\delta\,\tr(e_Ce_D)\big)^{2} f \big)\\
&=&\tr~ \eta\big(\big(\delta\,\tr(e_Ce_D)\big)^{2} f \big).
\end{IEEEeqnarray*}
As $\eta(\alpha f) = \eta(\alpha)f$ for $f$ a projection, and $\eta(\alpha^2)=2\alpha\eta(\alpha)$ for any scalar  $\alpha$, we get the following
\begin{IEEEeqnarray*}{lCl}
H\big(|\mathcal{F}(\Theta)|^2\big) 
&=&tr\big(\eta\big(\big(\delta\,\tr\big(e_Ce_D\big)\big)^{2}\big)f\big)\\
&=&\eta\big(\big(\delta\,\tr\big(e_Ce_D\big)\big)^{2}\big)\tr(f)\\
&=& 2 \delta \tr\big(e_Ce_D\big)\eta\big(\delta \hspace{0.1cm}\tr\big(e_Ce_D\big)\big)\tr(f)\\
&=& 2 \delta \tr\big(e_Ce_D\big)\eta\big(\delta\,\tr\big(e_Ce_D\big)\big) \frac{1}{[A:B]_0\tr(e_Ce_D)}\\
&=& \frac{2}{\delta}\eta\big(\delta\,\tr(e_Ce_D)\big),
\end{IEEEeqnarray*}
which completes the proof.\qed
\end{prf}

\begin{rmrk}\rm
It seems to be a difficult problem to  generalize the formula in \Cref{main} to the non-irreducible case.
\end{rmrk}

\begin{thm}\label{cocom}
Let $(B\subset C,D\subset A)$ be a quadruple of simple unital $C^*$-algebras with a conditional expectation from $A$ onto $B$ of index-finite type. If the quadruple is a co-commuting square, then 
\begin{center}
$H\big(|\mathcal{F}(\Theta)|^2\big)=\frac{2}{\delta}\eta\big(\frac{\delta}{[A:C]_0[A:D]_0}\big).$
\end{center}
\end{thm}
\begin{prf}
First observe that the quadruple $(A\subset C_1,D_1\subset A_1)$ is a commuting square. By \Cref{ppn4} and Proposition \ref{ppn1}, $q(C_1 ,D_1) =  \delta \hspace{0.1cm}e_{D_1} \star e_{C_{1}}$ is a projection. The following equations hold true:
\begin{IEEEeqnarray*}{lCl}
H\big(|\mathcal{F}(\Theta)|^2\big) &=& H\big(r^{2}  (e_{D_1}\star e_{C_1})^{2} \big)\\
&=& \tr\left(\eta\left(\frac{r^{2}}{\delta^{2}}\hspace{0.1 cm} q(C_1 , D_1)\right)\right)\\
&=& \tr\left(\eta\left(\frac{r^{2}}{\delta^{2}}\right)\,q(C_1 , D_1)\right)\\
&=& \eta\left(\frac{r^{2}}{\delta^{2}}\right)\tr\big(q(C_1 , D_1)\big)\\
&=& \frac{2r}{\delta}\eta\left(\frac{r}{\delta}\right)\frac{1}{r}\hspace{4 cm} (\mbox{by \Cref{fact1}})\\
&=& \frac{2}{\delta}\,\eta\left(\frac{r}{\delta}\right)\,. 
\end{IEEEeqnarray*}
Now, using the multiplicativity of the minimal index, it is easy to check that  $\displaystyle \frac{r}{\delta} = \displaystyle \frac{\delta}{[A:C]_0[A:D]_0}$, which concludes the proof.
%Therefore, we conclude that  $$H\big(|\mathcal{F}(\Theta)|^2\big)=\frac{2}{\delta} \hspace{0.1 cm} \eta\bigg(\frac{\delta}{[A:C]_0[A:D]_0}\bigg).$$
\qed
\end{prf}

\medskip
We conclude the paper with the following question.
\smallskip

\noindent \textbf{Problem:} Let $B\subset A$ be an irreducible inclusion of simple unital $C^*$-algebras with a conditional expectation from $A$ onto $B$ of index-finite type and suppose $C,D$ are two intermediate unital $C^*$-subalgebras. Can one provide a formula (in the spirit of \Cref{main}) for $H(\lvert \Theta\rvert ^2)?$

\section*{Acknowledgements}
The first author sincerely thanks Zhengwei Liu for useful exchange. K. C. Bakshi acknowledges support of INSPIRE Faculty grant DST/INSPIRE/04/2019/002754 and S. Guin acknowledges support of SERB grant MTR/2021/000818.

\bigskip

\bigskip

\noindent {\em Department of Mathematics and Statistics},\\
{\em Indian Institute of Technology Kanpur},\\
{\em Uttar Pradesh $208016$, India}
\medskip

\noindent {Keshab Chandra Bakshi:} keshab@iitk.ac.in, bakshi209@gmail.com\\
{Satyajit Guin:} sguin@iitk.ac.in\\
{Biplab Pal:} biplabpal32@gmail.com

\end{document}